\documentclass{article}
\usepackage{amssymb,amsmath}
\usepackage{graphics}

\def\qed{\hfill {\hbox{${\vcenter{\vbox{               
   \hrule height 0.4pt\hbox{\vrule width 0.4pt height 6pt
   \kern5pt\vrule width 0.4pt}\hrule height 0.4pt}}}$}}}

\def\bar{\overline}

\newtheorem{definition}{Definition}

\newtheorem{example}{Example}

{\qed\par\medskip}

\author{{\begin{tabular}{c} Conrad Creel \\
\small{\texttt{conradcreel@gmail.com}}\end{tabular}}
\and
{\begin{tabular}{c} Sam Nelson \\
\small{\texttt{knots@esotericka.org}}\end{tabular}}
\and
\small{\begin{tabular}{c}
Department of Mathematics, University of California, Riverside\\
900 University Avenue, Riverside, CA, 92521\end{tabular} }}

\date{}

\title{\Large \textbf{Symbolic computation with finite biquandles}}

\begin{document}
\maketitle

\begin{abstract}
A method of computing a basis for the second Yang-Baxter cohomology of a 
finite biquandle with coefficients in $\mathbb{Q}$ and $\mathbb{Z}_p$ from 
a matrix presentation of the finite biquandle is described. We
also describe a method for computing the Yang-Baxter cocycle
invariants of an oriented knot or link represented as a signed Gauss code.
We provide a URL for our \texttt{Maple} implementations of these algorithms.
\end{abstract}



\section{\large \textbf{Introduction}}

In \cite{HMN}, algorithms for computing the values of quandle counting 
invariants on virtual knots and for finding the finite quandles which 
define these invariants were described. In this paper we describe algorithms
for finding and computing a generalized form of this quandle counting invariant
using a similar methodology -- representing virtual knots and links 
as signed Gauss codes and representing finite biquandles as block matrices.

A \textit{biquandle} is a non-associative algebraic structure defined on a 
set $B$ whose axioms are motivated by thinking of the elements of $B$ as 
\textit{semi-arcs} in an oriented knot diagram and thinking of the four 
possible crossing relationships on inbound semi-arcs at a crossing as four 
binary operations on the set; the biquandle axioms are then the conditions
required in order to preserve the algebraic structure under the three 
Reidemeister moves. Biquandles have been studied in recent papers such as 
\cite{FJK}, \cite{KM}, \cite{NV} and more. The resulting algebraic structure 
is naturally a source of invariants of knots and links, both in the classical 
sense of disjoint unions of simple closed curves in $S^3$ and the 
combinatorial sense of Reidemeister equivalence classes of Gauss codes, also 
known as \textit{abstract knots} or \textit{virtual knots}
\cite{K2}, \cite{K1}.
 
Perhaps the simplest easily computable link invariant derivable from a finite 
biquandle is the counting invariant $|\mathrm{Hom}(K,T)|$ which counts 
homomorphisms from the knot biquandle $K$ into a finite target biquandle $T$. 
This counting invariant has been studied in the case of various finite 
biquandle structures defined algebraically (Alexander biquandles, 
quaternionic biquandles, etc.) in \cite{FJK} and in the case of biquandles
defined symbolically by biquandle matrices encoding the four operation tables
in \cite{NV}.

Simply counting homomorphisms, however, discards some information about the 
link
by treating all homomorphisms as the same. One way to incorporate information 
about the link diagram in the set of biquandle homomorphisms is to use the 
homomorphism to define a \textit{Boltzmann weight}, such as an integer 
power of a variable $t$, at each crossing; the product of these weights over 
all crossings defines an invariant of the colored knot diagram, and the sum of 
these weights over the set of all colorings defines a new knot invariant. The 
kinds of weight-assigning functions which make this construction work are 
elements of the second Yang-Baxter cohomology of the target biquandle as 
defined in \cite{C1}; hence, to find these invariants we must compute the 
second cohomology of our target biquandle.

In this paper we describe an algorithm for computing the second Yang-Baxter
cohomology of a finite biquandle with rational coefficients from its symbolic
matrix presentation. We give some examples of symbolic computation of the 
resulting Yang-Baxter cocycle invariants for virtual links represented by 
Gauss codes.

The paper is organized as follows: In section 2, we define biquandles and
describe ways of representing finite biquandles in \texttt{Maple}. In 
section 3 we describe a combinatorial generalization of knot theory known
as virtual knot theory and describe how virtual knots and links may be
represented in \texttt{Maple} as signed Gauss codes. In section 4 we 
describe our method for computing a basis for the second Yang-Baxter 
cohomology of a finite biquandle with coefficients in $\mathbb{Q}$ or 
a finite field $\mathbb{Z}_p$ and our method for computing the Yang-Baxter 
2-cocycle invariants of a knot or link using these cocycles. Our \texttt{Maple}
programs are available in the files \texttt{biquandles-maple} and
\texttt{yangbaxtercohomology.txt} downloadable at 
\texttt{www.esotericka.org/quandles}.

\section{\large \textbf{Finite biquandles}}

We begin this section with the definition of a biquandle \cite{KR}.

\begin{definition}
\textup{
A \textit{biquandle} is a set $B$ with four binary operations
$B\times B\to B$ denoted by
\[(a,b) \mapsto a^b, \ a^{\bar b}, \ a_b,\quad \mathrm{and} \quad a_{\bar b}\]
respectively, satisfying the following axioms:}
\begin{list}{}{}
\item[\textup{1.}]{\textup{For every pair of elements $a,b\in B$, we have}
\[
\mathrm{(i)} \ a=a^{b{\bar{b_a}}}, \quad
\mathrm{(ii)} \ b=b_{a{\bar{a^b}}}, \quad
\mathrm{(iii)} \ a=a^{\bar{b}b_{\bar a}}, \quad
\mathrm{and} \quad
\mathrm{(iv)} \ b=b_{\bar{a}a^{\bar b}}.
\]}
\item[\textup{2.}]{\textup{ Given elements $a,b\in B$, there are
elements $x,y\in B$, possibly but not necessarily distinct, such that}
\[
\mathrm{(i)} \ x=a^{b_{\bar x}}, \quad
\mathrm{ (ii) } \ a=x^{\bar b},\quad
\mathrm{ (iii)} \   b=b_{{\bar x}a}, \]\[
\mathrm{ (iv)} \  y=a^{\bar{b_y}},\quad
\mathrm{ (v) } \  a=y^b,   \quad \mathrm{and} \quad
\mathrm{ (vi)} \  b=b_{y{\bar a}}.
\]}
\item[\textup{3.}]{\textup{ For every triple $a,b,c \in B$ we have:}
\[
\mathrm{ (i)} \ a^{bc}=a^{c_bb^c}, \quad
\mathrm{(ii)} \ c_{ba} =c_{a^bb_a}, \quad
\mathrm{(iii)} \ (b_a)^{c_{a^b}}=(b^c)_{a^{c_b}}, \]
\[
\mathrm{(iv)} \ a^{{\bar b}{\bar c}}
=a^{\bar{c_{\bar b}}\bar{b^{\bar c}}},
\quad
\mathrm{(v)} \ c_{{\bar b}{\bar a}}
=c_{\bar{a^{\bar b}}{\bar {b_{\bar a}}}},
\quad \mathrm{and}  \quad
\mathrm{(vi)} \ (b_{\bar a})^{\bar{c_{{\bar {a^{\bar b}}}}}}
=(b^{\bar c})_{\bar{a^{\bar{c_{\bar b}}}}}.
\]}
\item[\textup{4.}]{\textup{Given an element $a\in B,$ there are elements
$x,y\in B$, possibly but not necessarily distinct, such that}
\[
\mathrm{ (i)} \   x=a_x, \quad
\mathrm{ (ii)} \  a=x^a, \quad
\mathrm{(iii)} \  y=a^{\bar y}, \quad \mathrm{and}  \quad
\mathrm{(iv) } \  a=y_{\bar a}.
\]
}
\end{list}
\end{definition}

If $B=\{x_1,x_2,\dots, x_n\}$ is a finite biquandle, then we define the
\textit{biquandle matrix} of $B$ to be the $2n\times 2n$ matrix
\[B=\left[
\begin{array}{r|r}
B_1 & B_2 \\ \hline
B_3 & B_4
\end{array}
\right]
\]
where the $i,j$ entry of the submatrix $B_m$ is $(B_m)_{i,j}=k$ where 
\[
(x_k)=\left\{
\begin{array}{ll}
(x_i)^{\bar{(x_j)}} & m=1\\
(x_i)^{(x_j)} & m=2 \\
(x_i)_{\bar{(x_j)}} & m=3\\
(x_i)_{(x_j)} & m=4. \\
\end{array}
\right.\]

To avoid confusing subscripts denoting elements of $B$ with the biquandle
operations, we generally drop the ``$x$''s and just write 
$B=\{1,2,\dots, n\}$.

In \texttt{Maple}, we represent a finite biquandle as a list of four 
$n\times n$ matrices $[M[1],M[2],M[3],M[4]]$ where $M[i]$ is the matrix
$B_i$ for $i=1,\dots, 4$. Thus, the biquandle word $1^{2_{\bar{3_2}}}$ in our
\texttt{Maple} notation is
\[M[2][1,M[4][2,M[3][3,2]]].\]

If axioms (1)-(3) are satisfied, $B$ is a \textit{birack}. The notation can be
simplified by denoting $S(a,b)=(b_a,a^b)$; then axiom (1) implies that
$S:B\times B\to B\times B$ is invertible with inverse given by $S^{-1}(a,b)=
(a^{\bar b},b_{\bar a})$. This is known as \textit{switch} notation;
a map $S$ is called a \textit{switch map} if it satisfies 
the \textit{set-theoretic Yang-Baxter equation}
\[(S\times \mathrm{Id})(\mathrm{Id}\times S)(S\times \mathrm{Id}) =
(\mathrm{Id}\times S)(S\times \mathrm{Id})(\mathrm{Id}\times S). \]
A biquandle is then a set $B$ with an invertible switch map 
$S:B\times B\to B\times B$ whose component maps additionally satisfy
axioms (2) and (4). See \cite{FJK} for more.

\begin{example}
\textup{The set $\mathbb{Z}_n=\{ 0,1,2,\dots, n-1\}$ is a biquandle with
operations} 
\[x^y=tz+(1-st)y \ \mod \ n, \quad x_y=sx \ \mod \ n,\] 
\[x^{\bar{y}}=t^{-1}z+(1-s^{-1}t^{-1})y \ \mod \ n, 
\quad x_{\bar{y}}=s^{-1}x \ \mod \ n\]
\textup{where $s,t$ are invertible elements of
$\mathbb{Z}_n$. This is an example of an \textit{Alexander biquandle}.
See \cite{KM} and \cite{LN} for more.}
\end{example}

In \cite{NV}, an algorithm is described for finding all biquandle structures 
on a set with a given finite cardinality, as well as an algorithm for counting
homomorphisms from a finitely presented biquandle described by a presentation 
matrix into the specified finite biquandle. The biquandle search algorithm 
works by taking a ``partially completed'' $2n\times 2n$ matrix with entries
in $\{0,1,2,\dots, n\}$ where a zero is considered a ``blank;'' the program
then selects a zero entry and replaces the zero with 
nonzero entries, propagating the value through the matrix using the biquandle 
axioms and equations obtained from them and searching for contradictions.
Zeroes are rated according to how many biquandle words in the list of axiom
testing conditions will be completed if the zero is filled in, via a program
called \texttt{ratezero}; the program selects a zero entry with maximal rating
in order to move through the search space more efficiently.
Any resulting matrices are appended to a working list, and the process 
repeats until all zeroes have been filled in. 

Our \texttt{Maple} code includes a program 
\texttt{abq} which finds the biquandle matrix for the Alexander biquandle 
$\mathbb{Z}_n$ with a choice of $s$ and $t$ invertible elements of 
$\mathbb{Z}_n$, as a way of generating biquandle matrices with larger 
cardinalities.

\section{\large \textbf{Virtual knots}}

Knots and links (disjoint unions of simple closed curves in $S^3$) are usually
represented as \textit{knot diagrams}, which are 4-valent graphs embedded in 
an oriented surface $\Sigma$ with vertices interpreted as crossings and 
decorated to indicate which strand passes over and which passes under.

Various methods exist for encoding the information from a knot diagram in a 
more code-friendly way. One such method uses \textit{Dowker-Thistlethwaite 
codes} or \textit{DT codes}, e.g. the 
\texttt{knotscape}\footnote{Available at 
\texttt{http://www.math.utk.edu/$\sim$morwen/knotscape.html}} package. For 
\textit{oriented links} in which each component has a preferred choice of 
direction, we use signed Gauss codes to represent link diagrams.

A \textit{signed Gauss code} for an oriented knot diagram $K$ is an ordered 
list of crossing labels (including over/under information) and signs 
encountered as one travels the knot the direction of the orientation. An 
oriented crossing is \textit{positive} if,
while looking in the positive direction of the overcrossing strand, the 
undercrossing strand is oriented right-to-left, and the crossing is 
\textit{negative} if the under-strand is oriented left-to-right. If we have 
an oriented link with multiple components, we separate the Gauss code 
components with commas.
For a given oriented link diagram, the corresponding signed Gauss code is 
well-defined up to ordering of the components and choice of base point for 
each component, with different choices of base point corresponding to cyclic 
permutations of the crossing labels. Given a signed
Gauss code obtained from a link, we can reconstruct the original link diagram 
up to local isotopy, i.e., isotopy of the surface $\Sigma$ in which the
link diagram is embedded.

To represent an $n$-crossing Gauss code in \texttt{Maple}, we use a vector 
whose components are Gaussian integers $\pm X$ for a positive crossing and 
$\pm(X+\sqrt{-1})$ for a negative crossing where $X\in\{1,\dots,n\}$, with 
the positive
entry representing the overcrossing label and the negative the undercrossing 
label. We use a ``0'' entry as an ``end of component'' indicator to separate
components in a link.

\begin{example}
\textup{The two-component link}
\raisebox{-0.4in}{\includegraphics{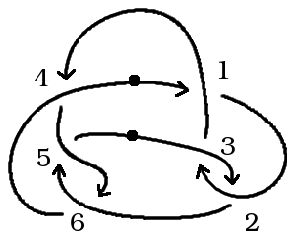}}
\textup{has signed Gauss code}
\[ [-1-I,2,-3,1+I,-4-I,5,-6,4+I,0,3,-2,6,-5,0]\]
\textup{in our \texttt{Maple} format.}
\end{example}

Knots and links may be defined combinatorially as equivalence classes of link 
diagrams under the equivalence relation generated by the Reidemeister moves. 
If we restrict our attention to Gauss codes in which every crossing number 
appears once as an overcrossing and once as an undercrossing, with the same 
sign for both instances, then equivalence classes of such codes under the 
equivalence relation generated by the Gauss code Reidemeister moves are 
called \textit{virtual knots}. These equivalence classes include some codes 
which do not correspond to knot diagrams which can be realized in the plane. 
Such Gauss codes may be interpreted as knot diagrams drawn on surfaces 
$\Sigma$ with genus, corresponding geometrically to 
knots or links in thickened surfaces $\Sigma\times [0,1]$ up to stabilization 
\cite{CKS}. When such knots are drawn on genus-0 paper, the crossings which
correspond to genus in the supporting surface, called \textit{virtual 
crossings}, are shown as circled intersections to distinguish them from the 
ordinary \textit{classical} crossings. In particular, virtual crossings have 
no over- or under-sense. Classical 
knots and links are then simply the subset of virtual knots and links which 
have representatives with supporting surface of genus 0.
See \cite{K1} for more about virtual knots.

For any oriented virtual link L there is an associated biquandle, the 
\textit{knot biquandle} $B(K)$, defined by assigning a generator for each 
\textit{semi-arc} (the portion of a link diagram going from one over or 
undercrossing point to the next) and two relations at each crossing:

\[\includegraphics{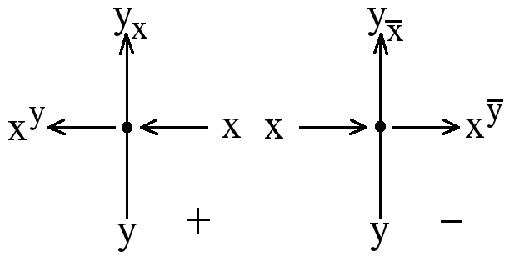}\]

Thus, for any virtual knot or link $K$, we have an associated biquandle 
$B(K)$ defined via a universal algebra-style presentation
with generators and relations; see \cite{KR} and \cite{FJK} for more. If we 
label semi-arcs with generator numbers $1, 2, \dots, n$ then such
a presentation can be expressed as a biquandle presentation matrix with 
non-zero entries corresponding to relations coming from the crossings.

A homomorphism from the knot biquandle $B(K)$ to a finite target biquandle $T$
can be regarded as a ``coloring'' of the knot diagram, i.e., an assignment of
an element of $T$ to every semiarc in the diagram such that above-pictured 
relations are satisfied at every crossing.

Initially, we used biquandle presentation matrices to represent knot
biquandles for the purpose of counting homomorphisms into finite biquandles.
The simplest method of computing $|\mathrm{Hom}(B(K),T)|$ is to generate a
presentation matrix for $B(K)$ from a signed Gauss code, then test all 
possible maps from the set $\{1, 2, \dots, 2n\}$ of generators of $B(K)$ to
the target biquandle $T=\{1, 2, \dots, m\}$ for the homomorphism
condition (i.e., that the images of the generators of the knot biquandle
satisfy the relations from the knot biquandle). However, for links with 
largish numbers of crossings this becomes very slow; an $n$-crossing link 
diagram has a knot biquandle with $2n$ generators, and brute-force checking 
all $(2n)^m$ maps where the target biquandle has cardinality $m$ is 
impractical. An improved method is to use a ``blank'' homomorphism and 
propagate values through using the homomorphism condition in a manner 
analogous to our biquandle finding program. This is implemented in our 
\texttt{Maple} code as \texttt{bhomlist}. 

An improved method of representing the knot biquandle of an oriented 
link is implemented in \texttt{bhomlist2} in the file
\texttt{yangbaxtercohomology.txt}, available on 
\texttt{www.esotericka.org/quandles}. Rather than using a biquandle 
presentation matrix, the program \texttt{gauss2list2} takes
a signed Gauss code and creates a list of biquandle relations using 
$M[1][i,j]$ for $i^j$, $M[2][i,j]$ for $i_j$, $M[3][i,j]$ for $i^{\bar{j}}$ 
and $M[4][i,j]$ for $i_{\bar{j}}$. The program \texttt{reducepreslist} then 
reduces the presentation by
looking for generators which appear on one side of the relation but not the 
other; such generators are then eliminated, with all instances being replaced 
by the equivalent word in the remaining generators. This results in a dramatic 
reduction in the amount of brute-force checking required to compute 
$|\mathrm{Hom}(B(K),T)|$. Unfortunately, this method does require brute-force 
checking of maps; we cannot fill in and propagate values through since all 
surviving generators appear on both sides of any relation in which they 
appear. Nevertheless, this method is computationally much faster than the 
method of \texttt{bhomlist} since most knot biquandle presentations can be 
reduced significantly.

\begin{example}
\[\includegraphics{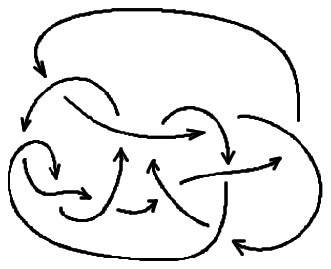}\]

\textup{
This 11-crossing diagram of the Conway knot has knot biquandle with 
22-generator presentation}
\begin{eqnarray*}
\langle 1,2,\dots, 22 & | & 1^{16}=2, \ 2_{15}=3, \ 3_8=4, \ 4^9=5,  \
5_{22}=6,\ 6^{\bar{11}}=7, \ \\ & & 
7_{\bar{20}}=8,\  8^3=9,\  9_4=10, \ 10^{\bar{21}}=11,\  
11_{\bar{6}}=12, \ \\ & & 
12_{\bar{17}}=13, \ 13^{\bar{18}}=14, \ 
14_{\bar{19}}=15, \ 15^2=16, \ 16_1=17, \\ & & 
17^{\bar{12}}=18, \ 18_{\bar{13}}=19, 
\ 19^{\bar{14}}=20, \ 20^{\bar{7}}=21, \\ & & 
\ 21_{\bar{10}}=22, \ 22^5=1\rangle.
\end{eqnarray*}

\textup{\texttt{gauss2pres2} reduces this to a 5-generator presentation with
generators $1, 8,$ $15, 16,$ and $21$ and relations too awkward to list here, 
though the interested reader can readily generate them with our code. 
For computing $|\mathrm{Hom}(K,T)|$ with $|T|=4$, this reduces the number of 
brute-force checks from $22^4$ to $5^4$, an improvement of 37,481\%; moreover, 
the number of relations which must be checked at each step is also reduced 
from 22 to 5.}
\end{example}

Given a homomorphism obtained from the reduced biquandle presentation from 
\texttt{gauss2pres2}, we can recover the corresponding coloring of the 
diagram by using \texttt{bhomcomplete}. We will need these complete colorings 
for computing the Yang-Baxter cocycle invariant in the next section.

\section{\large \textbf{Yang-Baxter cocycle invariants}}

For any link $L$ and finite biquandle $T$, the cardinality of the set
of biquandle homomorphisms from the knot biquandle $B(L)$ to $T$ is an
invariant of link type. A set is more than a mere cardinality, however, 
and we'd like to try to extract more information about the link type from 
the set of biquandles colorings.

One way to do this, described in \cite{C1}, is to use \textit{Yang-Baxter 
cocycles}. The authors define a homology theory for finite biquandles
with $n$-chains generated by ordered $n$-tuples of biquandle elements.
The boundary map is defined using a bijection between ordered $n$-tuples 
in $T^n$ and colorings of the $n$-cube graph $Q_n$ with biquandle elements 
defined by interpreting an $n$-tuple as a coloring of the
edges in a ``preferred path'' from the origin to $(1,1,\dots,1)$ 
in the $n$-cube graph $Q_n$. Such a coloring extends to a unique 
\textit{Yang-Baxter coloring} of the whole graph. The boundary map from 
standard cubical homology then induces via this bijection a boundary map
$\partial_n^{YB}:A[T^n] \to A[T^{n-1}]$ with coefficients 
in an abelian group $A$. The resulting homology and 
cohomology theories are the \textit{Yang-Baxter homology and 
cohomology} of the biquandle. We will primarily use $\mathbb{Q}$ as 
our coefficient ring for ease of computation;
our software includes both $\mathbb{Q}$ and $\mathbb{Z}_p$ versions of our 
programs for $p$ prime.

An alternative description of $H^2_{YB}(T;\mathbb{Q})$ in terms of knot 
diagrams makes it
clear how Yang-Baxter 2-cocycles can be used to retain some information from
the sets of biquandle colorings of a link diagram. Given a link diagram $L$
with a coloring by elements of $T$, we will define a function 
$\phi:T\times T \to \mathbb{Q}$ such that the sum of the values of $\phi$ on 
the inbound crossings at each positive crossing and on the outbound crossings 
at each negative crossing, always with the color on the undercrossing 
semi-arc listed first, is unchanged by the Reidemeister moves. The value
$\phi(x,y)$ at a crossing is called the \textit{Boltzmann weight} of the 
crossing. Invariance under
the second Reidemeister move is guaranteed by the convention just described;
invariance under the third Reidemeister move requires
\[\phi(x,y)+\phi(x^y,z)+\phi(y_x,z_{x^y}) =\phi(x,z_y)+\phi(y,z)
+\phi(x^{z_y},y^z)\]
which is precisely the condition that $\phi\in H^2_{YB}(T;\mathbb{Q}).$ 
\[\includegraphics{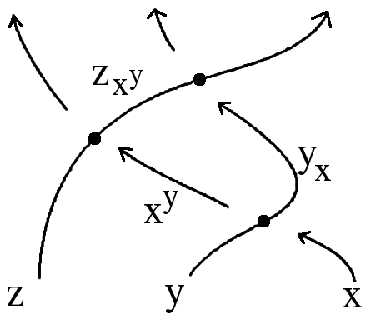}\quad 
\includegraphics{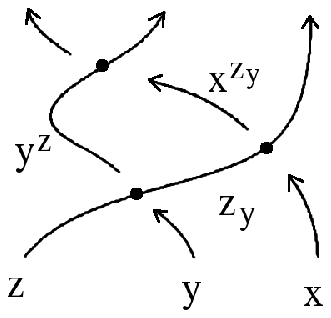}
\]

Invariance under the first Reidemeister move requires that $\phi(x,a)=0$ and
$\phi(a,y)=0$, where $x$ and $y$ are the elements associated to $a$ in 
biquandle axiom 4, for every $a\in T$. The subspace of 
$H^2_{YB}(T;\mathbb{Q})$ satisfying this condition is the 
\textit{reduced Yang-Baxter cohomology} of $T$ with coefficients in 
$\mathbb{Q}$. For each Yang-Baxter 2-cocycle 
$\phi:T\times T\to \mathbb{Q}$, the sum $\sum \phi(x,y)$ of 
Boltzmann weights over the set of all crossings in a biquandle-colored link 
diagram is unchanged by Reidemeister moves and can be understood as a kind of
``signature'' of the coloring homomorphism. 

For such a cocycle $\phi$, we compute the \textit{Yang-Baxter cocycle 
invariant} of a link $L$ by computing for each coloring 
$f\in \mathrm{Hom}(B(L),T)$ of $L$ by $T$
the sum over all crossings in $L$ of the Boltzmann weights $\phi(x,y)$ as 
above. Then, the set with multiplicities of these sums for all colorings of 
$L$ by $T$ is an invariant of links which contains more information than the 
cardinality of the set. For convenience, we can make each such sum the 
exponent
of a formal variable $t$ and then take the sum over the set of colorings;
there is no loss of information here since we are adding only powers of 
$t$ with coefficients of 1. Then,
\[\Phi_{YB}(L,T,\phi)=\sum_{f\in\mathrm{Hom}(B(L),T)} t^{\sum \phi(x,y)}.\]
Note that cohomologous cocycles define the same invariant.

To compute a basis for the rational cohomology vector space 
$H^2_{YB}(T;\mathbb{Q})$ from a biquandle matrix for $T$, we first make a 
matrix with rows corresponding to triples of biquandle elements and columns 
corresponding to characteristic functions 
$\chi_{(x,y)}(v,w)=\left\{\begin{array}{ll} 1 & (v,w)=(x,y) \\
0 & \mathrm{else}\end{array}\right.$ on pairs of biquandle 
elements; for each triple $(x,y,z)$ the entries in the columns corresponding 
to $\chi_{(x,y)}, \ \chi_{(x^y,z)}$ and $\chi_{(y_x,z_{x^y})}$ are set equal 
to 1, those corresponding to  $\chi_{(x,z_y)}, \ \chi_{(y,z)}$ and 
$\chi_{(x^{z_y},y^z)}$ are set equal to $-1$ and all other entries are $0$. 
$H^2_{YB}(T,\mathbb{Q})$ is then the null space of this matrix. To find a 
basis we use Gauss-Jordan elimination and then apply \texttt{getkernel} 
to read off a basis which is consistent between runs.

The program \texttt{ybcohom} then eliminates cohomologous cocycles by 
comparing the basis vectors from \texttt{ybcocycles} pairwise and keeping 
only one vector
from each pair whose difference is a coboundary. Finally, \texttt{redybcohom}
eliminates cocycles which do not satisfy the condition arising from the
type I Reidemeister move.

Our program \texttt{ybinv} takes as input a signed Gauss code, a biquandle 
matrix, and a Yang-Baxter 2-cocycle and computes the Yang-Baxter 2 cocycle 
invariant. Another program \texttt{ybinv2} takes a signed Gauss code and
biquandle matrix and computes a basis for the reduced cohomology of $T$, 
and for each such cocycle computes the resulting cocycle invariant, 
outputting a vector of cocycle invariant values. Our use of \texttt{getkernel}
ensures that the cocycles (and hence their invariants) are listed in the 
same order between runs.

Finally, we have modified versions of \texttt{ybinv} and associated programs, 
called \texttt{redybcohommodn}, \texttt{ybinvmodn}, 
etc. in which we replace 
$\mathbb{Q}$ with the finite field $\mathbb{Z}_p$ for $p$ prime.

\begin{example}
\textup{The trivial cocycle $\phi(x,y)=0$ has Yang-Baxter cocycle invariant 
$\Phi_{YB}(L,T,0)=|\mathrm{Hom}(B(L),T)|$ equal to the biquandle coloring 
invariant.}
\end{example}

\begin{example}
\textup{If the finite target biquandle T satisfies $a_b=a_{\bar b}=a$ for
all $a,b\in T$ then $T$ is a \textit{quandle} and the Yang-Baxter 2-cocycle
invariants associated to $T$ are the CJKLS state-sum invariants described
in \cite{CJKLS}.}
\end{example}

\begin{example}
\textup{This virtual knot is one of the \textit{Kishino knots} introduced
in \cite{KS}. The fact that it is not equivalent to the unknot, which was
shown in \cite{BF} using quaternionic biquandles, shows that 
the operation of connected sum for virtual knots, unlike the classical case, 
is dependent on the portion of the knot in which the sum is performed. The 
biquandle}
\[K=
\raisebox{-0.5in}{\includegraphics{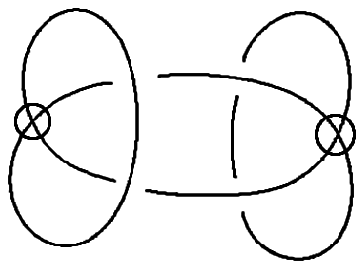}} \quad
T=\left[
\begin{array}{cccc|cccc}
1 & 4 & 2 & 3 & 1 & 3 & 4 & 2 \\
2 & 3 & 1 & 1 & 3 & 1 & 2 & 4 \\
4 & 1 & 3 & 2 & 2 & 4 & 3 & 1 \\
3 & 2 & 4 & 1 & 4 & 2 & 1 & 3 \\ \hline
1 & 3 & 4 & 2 & 1 & 4 & 2 & 3 \\
3 & 1 & 2 & 4 & 2 & 3 & 1 & 4 \\
2 & 4 & 3 & 1 & 4 & 1 & 3 & 2 \\
4 & 2 & 1 & 3 & 3 & 2 & 1 & 4 \\
\end{array}\right]\]
\textup{has reduced Yang-Baxter 2nd rational cohomology with basis 
$\{\phi_1,\phi_2\}$ where}
\[
\phi_1=-\chi_{(1,3)}-\chi_{(2,1)} -\chi_{(2,3)}-\chi_{(3,2)}, \]\[
\phi_2=\chi_{(1,3)}+\chi_{(1,4)}+\chi_{(2,1)}-\chi_{(2,3)}+\chi_{(3,1)}
-\chi_{(3,4)}.
\]
\textup{Our program \texttt{ybinv} computes the Yang-Baxter 2-cocycle 
invariant values 
\[\Phi_{YB}(K,T,\phi_1)=12+2t^{-1}+2t \quad \mathrm{and} \quad 
\Phi_{YB}(K,T,\phi_1)=12+2t^2+2t^{-2}.\] 
We note that both of these invariants distinguish the Kishino knot
$K$ from the unknot, which has a value of $\Phi_{YB}(K,T,\phi_i)=4$.}
\end{example}


\begin{thebibliography}{00}

\bibitem{BF}{S. Budden and R. Fenn.
The equation $[B,(A-1)(A,B)]=0$ and virtual knots and links.
\textit{Fundam. Math.} \textbf{184} (2004) 19-29.}

\bibitem{C1}{J. S. Carter, M. Elhamdadi, M. Saito. Homology theory for the 
set-theoretic Yang-Baxter equation and knot invariants from 
generalizations of quandles. \textit{Fund. Math.} \textbf{184} (2004) 31-54.}

\bibitem{CJKLS}{J. S. Carter, D. Jelsovsky, S. Kamada, L Langford and M. Saito.
Quandle cohomology and state-sum invariants of knotted curves and surfaces. 
\textit{Trans. Am. Math. Soc.} \textbf{355} (2003) 3947-3989.}

\bibitem{CKS}{J. S. Carter, S. Kamada, M. Saito. Stable Equivalence of Knots 
on Surfaces and Virtual Knot Cobordisms, \textit{J. Knot Theory Ramifications} 
\textbf{11} (2002) 311-322.}

\bibitem{DT}{C. H. Dowker and M. B.  Thistlethwaite.
On the classification of knots.
\textit{C. R. Math. Acad. Sci., Soc. R. Can.} \textbf{4} (1982) 129-131.}

\bibitem{FJK}{R. Fenn, M. Jordan-Santana, and L. Kauffman. 
Biquandles and virtual links.  \textit{Topology Appl.} \textbf{145} (2004) 
157-175.}

\bibitem{HMN}{R. Henderson, T. Macedo and S. Nelson. Symbolic Computation 
with finite quandles. \textit{J. Symbolic Comput.} \textbf{41} (2006) 811-817.}

\bibitem{K1}{L. H. Kauffman. Virtual knot theory. \textit{European J. Combin.}
\textbf{20}  (1999) 663-690.}

\bibitem{K2}{N. Kamada and S. Kamada. Abstract link diagrams and virtual knots.
\textit{J. Knot Theory Ramifications}  \textbf{9}  (2000) 93-106.}

\bibitem{KM}{L. H. Kauffman and V. O. Manturov. Virtual biquandles.  
\textit{Fund. Math.}  \textbf{188}  (2005) 103-146.}

\bibitem{KR}{L. H. Kauffman and D. Radford. Bi-oriented quantum algebras, 
and a generalized Alexander polynomial for virtual links. \textit{Contemp. 
Math.} \textbf{318} (2003) 113-140.}

\bibitem{KS} T. Kishino and S. Satoh. A note on non-classical virtual knots.
\textit{J. Knot Theory Ramifications} \textbf{13} (2004) 845-856.

\bibitem{LN}{D. Lam and S. Nelson. Classsification of finite Alexander 
biquandles. Preprint available at arXiv.org:math.QA/0611887}

\bibitem{NV}{ S. Nelson and J. Vo. Matrices and finite biquandles. 
\textit{Homology, Homotopy and Applications} \textbf{8} (2006) 51-73.
}

\end{thebibliography}
\end{document}